\documentstyle[amssymb]{article}
\title{\textsf{The Andrews-Curtis Conjecture\\ and\\ Black Box Groups}}

\author{\textsf{ Alexandre V. Borovik}\thanks{Supported by the Royal Society and
The Leverhulme Trust.}  \and \textsf{Evgenii I. 
Khukhro}\thanks{Supported by the Russian Foundation for Basic 
Research, grant no. 99-01-00576, and by a grant of Russian 
Ministry for Education in fundamental  research no.\ E00-1.0-77. 
} \and \textsf{ Alexei G. Myasnikov\thanks{Partially supported 
by  LMS Scheme 4 grant 4523 and by NSF grant DMS-9973233}}}
\date{\textsf{31 October 2001}}

\newtheorem{theorem}{Theorem}[section]
\newtheorem{lemma}[theorem]{Lemma}
\newtheorem{fact}[theorem]{Fact}
\newtheorem{corollary}[theorem]{Corollary}
\newtheorem{conjecture}{Conjecture}
\newtheorem{question}[conjecture]{Question}
\newtheorem{proposition}[theorem]{Proposition}

\begin{document}

\pagestyle{myheadings} \markright{{\scriptsize \textsf{A. V. Borovik,
 E. I. Khukhro and A. G. Myasnikov} $\bullet$ \textsf{ The Andrews-Curtis
Graph} $\bullet$ 31.10.01}}

\maketitle

\begin{abstract}
The paper discusses the Andrews-Curtis graph $\Delta_k(G,N)$ of a
normal subgroup $N$ in a group $G$. The vertices of the graph are
$k$-tuples of elements in $N$ which generate $N$ as a normal
subgroup; two vertices are connected if one them can be obtained
from another by certain elementary transformations. This object
appears naturally in the theory of black box finite groups and in
the Andrews-Curtis conjecture in algebraic topology \cite{AC}.

We suggest an approach to the Andrews-Curtis conjecture based on
the study of  Andrews-Curtis graphs of finite groups, discuss
properties of Andrews-Curtis graphs of some classes of finite
groups and results of computer experiments with generation of
random elements of finite groups by random walks on their
Andrews-Curtis graphs.

\end{abstract}

\section{Introduction}

The concept of a {\it black box group}  is a formalisation of a
probabilisitc approach to computational problems of finite group
theory.  For example, given two square matrices $x$ and $y$ of
size, say, 100 by 100 over a finite field, it is unrealistic to
list all elements in the group $X$ generated by $x$ and $y$ and
determine the isomorphism class of $X$ by inspection.  But this
can often be done, with an arbitrarily small probability of error,
by treating $X$ as a black box group, that is, by studying a
sample of random products of the generators $x$ and $y$.  The
explosive growth of the theory of black box groups in recent years
is reflected in numerous publications (see, for example, the
survey paper \cite{L-G} on the computational matrix group project)
and the fundamental work \cite{kantor}), and algorithms
implemented in the software packages GAP \cite{G} and MAGMA
\cite{BCP}.  A  critical  discussion of the concept of black box
group can be found in \cite{babai2}, while \cite{babai-beals}
contains a detailed survey of the subject.

In this paper, we look at the problem of generating random
elements from a normal subgroup of the black box groups. The
underlying structure, a version of the product replacement graph,
has rather unexpectedly happened to be the Andrews--Curtis graph
which appears in a certain group-theoretic problem of algebraic
topology.

We set the scene in Section~\ref{sec:PR} where we
briefly survey the known results about the product replacement
algorithm, and in Section~\ref{subsec:AC}, where we introduce the
Andrews--Curtis graph. If $G$ is a group (not necessary finite)
and  $N \lhd G$, then the {\em Andrews--Curtis graph\/}
$\Delta_k(G,N)$  is the graph whose  vertices are $k$-tuples of
elements in $N$ which generate $N$ as a normal subgroup:
$$ \left\{\, (h_1,\ldots, h_k) \mid \langle
h_1^G,\ldots, h_k^G\rangle = N \,\right\}.$$  Two vertices are
connected by an edge if one of them is obtained from another by
one of the moves:
\begin{eqnarray*}
(x_1, \ldots, x_i, \ldots, x_j, \ldots, x_k) & \longrightarrow & (x_1, \ldots,x_ix_j^{\pm 1}, \ldots,  x_k), \; i \ne j \\[2ex]
(x_1, \ldots, x_i, \ldots, x_j, \ldots, x_k) & \longrightarrow & (x_1, \ldots,x_j^{\pm 1}x_i, \ldots,  x_k), \; i \ne j \\[2ex]
(x_1, \ldots, x_i, \ldots, x_k) & \longrightarrow &  (x_1,\ldots,x_i^{w}, \ldots, x_k), \;  w\in G\\[2ex]
(x_1, \ldots, x_i, \ldots, x_k) & \longrightarrow & (x_1,
\ldots,x_i^{-1}, \ldots,  x_k).
\end{eqnarray*}

The well-known Andrews--Curtis Conjecture provides the main source
of motivation for the paper:
\begin{quote}
\textsf{\sl For a free group $F_k$ of rank \/ $k \geqslant 2$, the
Andrews--Curtis graph\/ $\Delta_k(F_k, F_k)$ is connected.}
\end{quote}
The main goal of the paper is to suggest a possible approach to
construction of a counterexample to the Andrews--Curtis conjecture
using finite groups (Section~\ref{sec:attack}). To that end, we
need a good understanding of the Andrews--Curtis graphs of finite
groups.
 In Section~\ref{sec:3}, we derive bounds for the diameters of the
Andrews-Curtis graphs $\Delta_k(G,G)$ of finite simple groups $G$.
Section~\ref{sec:soluble} contains a detailed study of the
Andrews-Curtis graphs of finite soluble groups. In particular, we
give a fairly complete description of the connected components  of
the graphs $\Delta_k(G,G)$ for finite nilpotent and soluble
groups $G$ (Theorems~\ref{th:nilpotent} and \ref{th:soluble}).

In Section~\ref{sec:ACWalks} we  return to the black box group
setting and use random walks on $\Delta_k(G,N)$ for generating
pseudorandom elements of a normal subgroup $N$ of a black box group $G$
and discuss the practical performance of our algorithm.
The Andrews-Curtis
graph has the following apparent advantages over the commonly used product
replacements graph $\Gamma_k(N)$ (it described in Section~\ref{sec:PR}).
The fact that the diameter
${\rm diam}\,(\Delta_k(G,N))$ is much smaller than that of
$\Gamma_k(N)$ suggests the possibility that the mixing time of a
random walk on $\Delta_k(G,N)$ is smaller than the mixing time of
a random walk on $\Gamma_k(N)$. The vertices of $\Delta_k(G,N)$
are all tuples in $N^k \smallsetminus (1,\ldots, 1)$, hence the
sample of elements of $N$ obtained by taking random components of
random generating tuples (vertices) in $\Delta_k(G,N)$  is not
biased.

\section{The Product Replacement Algorithm}
\label{sec:PR}

\subsection{A brief survey}

 A problem which we immediately encounter when
dealing with black box groups is how to construct a good black box
for the subgroup generated by given elements. For example, given a
group generated by a collection of matrices,
\begin{eqnarray*}
{ X} & \leqslant & GL_N({\mathbb F}_q),\\[2ex]
{ X} & = & \langle x_1, \ldots, x_k \rangle
\end{eqnarray*}
how can we produce (almost) uniformly distributed independent
random elements from ${X}$? The commonly used solution is
the {\em product replacement algorithm} (PRA) \cite{celler2}.

Denote by $\Gamma_k({X})$ the graph whose vertices are
generating $k$-tuples of elements in ${X}$ and edges are
given by the following \textit{Nielsen transformations}
\cite{mks}:
\begin{eqnarray*}
(x_1, \ldots, x_i, \ldots, x_k)  & \longrightarrow & (x_1,
\ldots, x_j^{\pm 1}x_i, \ldots, x_k)\\[1ex]
(x_1, \ldots, x_i, \ldots, x_k)  & \longrightarrow & (x_1, \ldots,
x_ix_j^{\pm 1}, \ldots, x_k),
\end{eqnarray*}
where $j \neq j$. Sometimes it is more convenient to consider the
extended graph ${\tilde \Gamma}_k({X})$ which has extra
edges corresponding to the transformations
\begin{eqnarray*}
(x_1, \ldots, x_i, \ldots, x_k)  & \longrightarrow & (x_1, \ldots,
x_i^{-1}, \ldots, x_k).
\end{eqnarray*}
The recipe for production of random elements from $ X$ is
deceptively simple: walk randomly over this graph and select
random components $x_i$. The detailed discussion of theoretical
aspects of this algorithm can be found in Igor Pak's survey
\cite{P2}. Pak \cite{P3} has also shown  that, if $k$ is
sufficiently big, the mixing time for a  random walk on
$\Gamma_k({X})$ is polynomial in $k$ and $\log
|{X}|$. Here, {\em mixing time} $t_{\rm mix}$ for a random
walk on a graph $\Gamma$ is the minimal number of steps such that
after these steps
$$
\frac{1}{2} \sum_{v\in \Gamma} \left| P(\hbox{get at } v) -
\frac{1}{\#\Gamma}\right| < \frac{1}{e}.
$$
At the intuitive level, this means that the distribution of the
end points of random walks on $\Gamma$ is sufficiently close to
the uniform distribution.

The graph $\Gamma_k({X})$ is still a very mysterious
object. Notice, in particular, that, in general, the graph
$\Gamma_k(G)$ is not connected. However, an elementary argument,
 due to Babai,
shows that $\Gamma_k(G)$ is connected for $k> 2 \log_2 |G|$. In
that case  the diameter $d(\Gamma_k(G))$ can be bounded as
$$d(\Gamma_k(G)) < C\cdot \log^2 |G|.$$

The following very natural question is still open.

\begin{conjecture} \textsf{ If\/ $G$ is a finite simple group,
the graph\/
$\Gamma_k(G)$ is connected for\/ $k \geqslant 3$.}
\end{conjecture}

Observe, that if we denote by $d(G)$ the minimal number of
generators for $G$, then for {\it every  finitely generated group}
 $G$ and every $k \geq d(G)+ 1$ the graph ${\tilde \Gamma}_k(G)$  is
connected if and only if the graph $\Gamma_k(G)$ is connected
(\cite{P2}).

A conceptual explanation of the good properties of the product
replacement algorithm is provided by the remarkable observation by
Lubotzky  and Pak.

\begin{theorem} {\rm (Lubotzky  and Pak \cite{LuP})}
If\/ ${\rm Aut} \, F_k$ satisfies Kazhdan property {\rm(}T{\rm )},
then mixing time $t_{\rm mix}$ of a random walk on a
component of\/ $\Gamma_k(G)$ is bounded as
$$
t_{\rm mix} \leqslant C(k) \cdot \log_2 |G|.
$$
\label{th:LuP}
\end{theorem}

Thus the issue is reduced to the long standing conjecture:

\begin{conjecture}
\textsf{ For\/ $k \geqslant 4$,  ${\rm Aut} \, F_k$ has} {\rm
(}T{\rm )}.
\end{conjecture}

Following \cite{K},
we say that a topological group $G$ satisfies the Kazhdan (T)-property if, for some
compact set  $Q \subset G$,
$$
K = \inf_{\rho} \inf_{v \ne 0} \max_{q \in Q}
\frac{\|\rho(q)(v)-v\|}{\|v\|} > 0,
$$
where $\rho$ runs over all unitary representations of $G$ without
fixed non-zero vectors. In our context, ${\rm Aut} \, F_k$ is endowed with the discreet topology.

Summarizing our brief discussion of PRA it worthwhile to mention
here that despite on computer experiments which show a good
overall  performance of PRA  still there are two major theoretical
obstacles when running PRA:

(O1) Connectivity of the graph $\Gamma_k(G)$;

(O2) The bias in the output of PRA.

It seems, the both obstacles can be removed by taking  $k$ large
enough \cite{P2}.  But then this increases the size of the
generating set thus affecting the performance of PRA.

\subsection{Normal subgroups of black box groups}
\label{sec:norblack}

 Assume that we know that some elements $y_1,\ldots,
y_k$ of a black box group ${ X}$ belong to a proper normal
subgroup of ${ X}$.

\begin{question}
\textsf{ How one can construct a good black box for the normal
closure
$${ Y} = \left\langle y_1^{ X},\ldots, y_k^{ X}
\right\rangle?$$ }
\end{question}
One possibility is to run a random walk on the Cayley graph for
$ Y$ with respect to the union of the conjugacy classes
$$y_1^{ X} \cup \cdots\cup y_k^{ X}$$
as the generating set for $ Y$.

Notice, that if we know that $ Y$ is a simple group then a
deep result by Liebeck and Shalev \cite[Corollary~1.14]{LiSh}
asserts that, for a finite simple group $G$ and a conjugacy class
$C \subset G$, the mixing time of the random walk on the Cayley
graph ${\mathsf Cayley}(G,C)$  is at most $c\log^3|G|/\log^2|C|$.

There are two remarks in order here:

1) As numerous  experiments showed, in general,  PRA performs much
better then any standard random walk on a Cayley graph of the
subgroup $N$ (\cite{P2});

2) Even though the mixing time of the PRA is polynomial (in
cardinality of the generating set and $\log_2|Y|$),
it is not a priori clear  how
many random elements $y_i^x$ one has to take to form a generating set of $Y$.

 These observations (not to mention the general obstacles (O1) and (O2))
 encourage one to look for other methods for constructing black box generator
for normal subgroups of black box groups. In the next section we
wish to discuss a modification of a product replacement algorithm
whose practical performance as a black box generator for normal
subgroups is better than a random walk or PRA with respect to the
generating set $ y_1^{ X} \cup \cdots\cup y_k^{ X}.$

\section{Andrews--Curtis graphs}
\label{subsec:AC}

Let $G$ be a group and $N \lhd G$.  Denote by $V_k(G,N)$ the set
of all $k$-tuples of elements in $N$ which generate $N$ as a
normal subgroup of $G$:
$$V_k(G,N) =  \left\{\, (h_1,\ldots, h_k) \in G^k \mid \langle
h_1^G,\ldots, h_k^G\rangle = N \,\right\}.$$ Of course, if the
group $N$ is simple then  $V_k(G,N) = N^k \smallsetminus
\{(1,\ldots, 1)\}$.

We define the {\em Andrews--Curtis graph\/} $\Delta_k(G,N)$ as the
graph with the set of vertices $V_k(G,N)$  and such that two
vertices are connected by an edge if one of them is obtained from
another by one of the following moves (\textit{elementary
Andrews--Curtis transformations}, or \textit{AC-transformations}):
\begin{eqnarray*}
(x_1, \ldots,x_i, \ldots,  x_k) & \longrightarrow & (x_1,
\ldots,x_ix_j^{\pm 1}, \ldots,
x_k), \; i \ne j, \\[2ex] (x_1, \ldots, x_i, \ldots,  x_k) & \longrightarrow & (x_1,
\ldots,x_j^{\pm 1}x_i, \ldots,  x_k), \; i \ne j, \\[2ex] (x_1, \ldots, x_i, \ldots,  x_k) &
\longrightarrow &  (x_1, \ldots,(x_i^{-1}), \ldots, x_k),
\\[2ex]
(x_1, \ldots, x_i, \ldots,  x_k) & \longrightarrow &  (x_1,
\ldots,x_i^w, \ldots,  x_k),
 \; w\in G. \end{eqnarray*}

Notice that the moves are invertible and thus give rise to a
non-oriented graph.

Sometimes it is convenient to consider a modification of the graph
$\Delta_k(G,N)$. Namely, if $A$ is a given finite set of
generators for $G$ then the graph $\Delta_k(G,N,A)^*$ has the same
set of vertices $V_k(G,N)$ which are connected by the same edges
as above, provided only that $w \in A$. In this case,  the number
of edges adjacent to a given vertex is finite (even if the group
$G$ is infinite).

Observe also, that if the group $G$ is abelian then $\Delta_k(G,N)
= {\tilde \Gamma}_k(N)$. Moreover, if $Ab(G)$ is  abelianization
of $G$, i.e., $Ab(G) = G/[G,G]$, then the canonical epimorphism $G
\rightarrow Ab(G)$ induces an adjacency-preserving map of graphs
$$\Delta_k(G,N) \rightarrow  {\tilde \Gamma}_k(Ab(G)).$$

The name and initial motivation to study graphs $\Delta_k(G,G)$
comes from the Andrews-Curtis Conjecture (1965) (AC-conjecture):

\begin{conjecture} {\rm ( Andrews and Curtis \cite{AC})}
\textsf{For $k \geqslant 2$, the Andrews--Curtis graph
$\Delta_k(F_k, F_k)$ is connected.}
\end{conjecture}

Obviously, for every group $G$ the graph $\Delta_k(G, G)$ is
connected if and only if the graph $\Delta_k(G, G)^*$ is
connected.

There is an extensive literature on the subject, see for example,
\cite{AK,BuM,HM}. Still virtually nothing is known about the
properties of the Andrews--Curtis graph for free groups. Some
potential counterexamples to the AC-conjecture (originated in
group theory and topology) were recently  killed by application of
{\em genetic algorithms} \cite{My1}, \cite{My2}. But the most
formidable stand untouched.

 One of the possible approaches to the AC-conjecture is based on the study
 of Andrews-Curtis graphs of  quotients of $F_k$ which are "close" to $F_k$.
 This is one of the few known positive results on connectivity of
 Andrews-Curtis graphs of relatively-free groups:

\begin{fact} {\rm (A.~G.~Myasnikov \cite{My})}
For the free soluble group $F_n^{(m)}$ of class $m$ and all\/ $k
\geqslant n$, the Andrews--Curtis graph\/ $\Delta_k(F_n^{(m)},
F_n^{(m)})$ is connected.  \label{fact:myasnikov}
\end{fact}

In Section \ref{sec:attack} we suggest a possible line of attack
at this problem which involves the study of the Andrews--Curtis
graphs $\Delta_k(G,G)$ for finite groups $G$.

\section{Random walks on Andrews--Curtis graphs }
\label{sec:ACWalks}

\subsection{AC-replacement algorithm}

In this section we discuss random walks on $\Delta_k(G,N)$ as an
alternative approach to black box generators of elements from a
normal subgroup $N$ of a group $G$.

Let $G$ be  a finite group and $N \lhd G$. If the graph
$\Delta_k(G,N)$  is connected then a nearest neighbour random walk
on this graph is an irreducible aperiodic Markov chain. Hence by
Perron-Frobenius theory it has uniform equilibrium distribution.

This suggest the following modification of the PRA which we call
AC-replacement algorithm ($ACR_k(G,N)$): run a nearest neighbour
random walk on $\Delta_k(G,N)$ for $t$ steps and return a random
component of the tuple in the stopping state.

\begin{conjecture} \textsf{ Let $G$ be a black box group and $N \lhd G$.
The AC-replacement algorithm $ACR_k(G,N)$ provides  a `good' black
box for $N$ (at least for some $k$).}
\end{conjecture}

In practice, a modification of the process,
when the the last changed component of the generating tuple (say, $x_ix_j^{\pm 1}$)  is
multiplied into the cumulative product $x$,
appears to be more effective:

\begin{itemize}
\item Initialise $x :=1$.
\item Repeat
    \begin{itemize}
        \item[$\circ$] Select random $i \ne j$ in $\{\,1,\ldots, k\,\}$.
        \item[$\circ$]
            \begin{itemize}
                \item With equal probabilities, replace $x_i := x_ix_j^{\pm 1}$ or $x_i := x_j^{\pm 1}x_i$, or
                \item produce random $w\in G$ and replace
                    $$x_i := x_i(x_j^w)^{\pm 1} \hbox{ or } x_i := (x_j^w)^{\pm 1}x_i.$$
             \end{itemize}
        \item[$\circ$] Multiply $x_i$ into $x$: $$x := x\cdot x_i.$$
    \end{itemize}

\item Use $x$ as the running output of a black box for $N$.
\end{itemize}

Using results on Markov chains, Leedham-Green and O'Brien
\cite{L-GO'B} had shown that  the distribution of the values of
the cumulative product $A$ converges exponentially to the uniform
distribution on $N$. However, the issue of explicit estimates is
open and represents a formidable problem.

In the next section we report on some computer experiments which
support the conjecture above.

\subsection{Generation of random elements\\ in simple normal
subgroups:\\ computer experiments} \label{sec:exp}

Here we give a brief discussion of some  computer experiments
related to the normal subgroups of black box  groups.

We run only a limited number of experiments, concentrating on the
generation of the alternating group ${\rm Alt}_n$ as a normal
subgroup of ${\rm Sym}_n$ by very short elements, for example, by
the involution $(12)(34)$ or by a $3$-cycle $(123)$. The two
series of  experiments were run, correspondingly, by the first
author in GAP \cite{G} and by Alexei D. Myasnikov (City College,
New York) using bespoke C++ code. We looked at the distribution of
numbers of cycles in the random permutation produced by
\begin{itemize}
\item the random walk on the Andrews-Curtis graph $\Delta_k({\rm
Sym}_n, {\rm Alt}_n)$,

\item  by the standard product replacement
algorithm,

\item and by a random walk on the Cayley graph ${\mathsf Cayley}({\rm
Alt}_n, x^{{\rm Sym}_n})$.
\end{itemize}

This particular criterion was chosen because of the importance of
permutations with small number of cycles in black box recognition
algorithms for symmetric groups \cite{BrP}.

We used in our experiments  the AC-replacement  algorithm (ACR)
with and without the cumulative product, as described in
Section~\ref{subsec:AC}.

Our experiments with alternating groups of degrees varying from
$10$ to $100$ have shown that when the generator is ``small'', a
very good convergence of the sample distribution to the uniform
distribution was achieved after $k\cdot n \cdot [\log_2 n]$ steps
of the algorithm, even if we worked with a very short generating
tuple, $k= 2$, $3$ or $5$.

The degree of convergence was measured by comparing the
distribution of the numbers of cycles in the cycle decomposition
of a permutation $x$ produced by the AC-replacement algorithm with
the theoretical distribution (easily computable from the Stirling
numbers of the first kind), and also by comparing the distribution
of the values $1^x$ with the uniform distribution on the set
$\{1,\ldots,n\}$. In both cases we used the $\chi^2$ criterion
with the significance level 95\%.

The performance of the ACR algorithm was, as a rule, better than a
random walk on the Cayley graph with respect to a conjugacy class
of the generator.

Also, the use of the cumulative product significantly improved
performance of the algorithm.

The standard product replacement algorithm has shown a very good
performance when the generating tuple was sufficiently long, or
when the initial generating tuple $(g_1,\ldots, g_k)$ was chosen
at random.

However, the ACR algorithm has shown robustness with respect to
choice of very small generators. This property is valuable in
certain applications, when one should expect to deal with the
initial generating tuple which is not representative of the
`average' elements of the normal subgroup. For example,
computations with centralizers of involutions of the type done in
\cite{altseimer-borovik,borovik} require computation of the normal
subgroups generated by involutions.

A discussion of similar experiments can be found in \cite{bl-gnf}.

\section{\relax Andrews--Curtis graphs\\ of finite groups}
\label{sec:finite}

\subsection{General bounds}
\label{sec:general}

 We show here that  if $k$ is large enough then the AC-graph $\Delta_k(G,G)$
of a finite group $G$ is connected. The proof is easy and similar
to the analogous result for the graph $\Gamma_k(G)$ (though
estimates are better).

Let $nd(G)$ be the minimal number of elements needed to generate
$G$ as a normal subgroup. Let $nd_m(G)$ be the maximal size of a
minimal set of normal generators of $G$ ($Y$ is a minimal set of
normal generators for $G$ if $\langle Y^G \rangle = G$, but
$\langle Y_0^G \rangle \neq  G$ for every proper subset $Y_0$ of
$Y$).

\begin{proposition}
\label{pr:confin} Let\/ $G$ be a finite group $G$. If\/ $k \geq nd(G)
+ nd_m(G)$ then the graph $\Delta_k(G,G)$ is connected.
\end{proposition}

{\it Proof.} Let $n = nd(G)$, $n_m = nd_m(G)$, and $k \geq n+n_m$.
Denote by $V_t(G)$ the set of all $t$-tuples which generate $G$ as
a normal subgroup.  Fix a tuple $h = (h_1, \ldots, h_n) \in
V_n(G)$.   Then $k$-tuple $h(k) = (h_1, \ldots, h_n, 1, \ldots,1)$
 is in $V_k(G)$. Now if $g = (g_1, \ldots, g_k) \in V_k(G)$, then
 there are  $n_m$ components of $g$, say $g_{n+1}, \ldots, g_k$, such that
 $(g_{n+1}, \ldots, g_k) \in V_{n_m}(G)$. It follows, that $g$ is connected
 in $\Delta_k(G,G)$ to $(h_1, \ldots,h_n, g_{n+1}, \ldots, g_k)$. Obviously,
 the latter one is connected to $h(k)$. Hence any tuple $g \in V_k(G)$ is
 connected in $\Delta_k(G,G)$ to $h(k)$, so the whole graph is
 connected. \hfill $\Box$

\subsection{\relax Andrews--Curtis graphs\\ of finite simple
groups} \label{sec:3}

In this section we give good and easy estimates (modulo known hard
results) of diameters of AC-graphs of finite simple groups.

\begin{theorem} If\/ $G$ is a finite simple
group and\/ $k \geqslant 2$ then the graph $\Delta_k(G,G)$ is
connected and
$$
{\rm diam}\,(\Delta_k(G,G)) < c\cdot k \cdot \log |G|
$$
for some constant\/ $c$.
\end{theorem}

This is a very crude estimate; the proof of the theorem contains
many possible directions for improvement, see
Proposition~\ref{prop:refined} below.

\paragraph{Proof.} If $G$ is a finite simple group,
then the {\em covering number\/} ${\rm cn}(G)$ is defined as
$$
{\rm cn}(G) = \min \{\,n \mid C^n = G \hbox{ for every conjugacy
class } C \subset G\,\}.
$$
The {\em Ore constant} ${\rm or}(G)$ is defined as
$$
{\rm or}(G) = \min \{\,n \mid C^n = G \hbox{ for some conjugacy
class } C \subset G\,\}.
$$
The prominent Ore-Thompson Conjecture asserts that ${\rm or}(G) =
2$ for all finite simple groups $G$.

The theorem follows from a simpler proposition.

\begin{proposition} If\/ $G$ is a finite simple group
and\/ $k \geqslant 2$ then
$\Delta_k(G,G)$ is connected  and
$$
{\rm diam}\,(\Delta(G,G)) \leqslant 4(k\cdot {\rm or}(G) + {\rm
cn}(G)).
$$
\label{prop:refined}
\end{proposition}

\paragraph{Proof of Proposition~{\ref{prop:refined}}.}
Set $d = {\rm or}(G)$. Let ${\mathsf{x}} = (x_1,1,\ldots, 1)$ be a
vertex in $\Delta_k(G,G)$ with $x_1$ chosen from the conjugacy
class $C$ such that $C^d = G$. We shall prove that $\mathsf{x}$
can be connected to an arbitrary vertex ${\mathsf{y}} =
(y_1,\ldots, y_k)$  of $\Delta_k(G,G)$.

Indeed, since $y_i = x_1^{w_1}\cdots x_1^{w_d}$ for some $w_j \in
G$, we get, after application to the tuple $(x_1,\ldots, x_k)$ of
$d(k-1)$ pairs of moves of the form (below $w_0 = 1$):
\begin{eqnarray*}
(z_1, \ldots, z_k) & \longrightarrow &  (z_1^{w_{i-1}^{-1}w_i},
z_2, \ldots,
z_k), \; i = 1,2,\ldots, d \\[1ex]
(z_1, \ldots, z_k)  & \longrightarrow & (z_1, \ldots, z_iz_1,
\ldots,  z_k)
\end{eqnarray*}
the tuple $(x_1^{w},y_2, \ldots,y_k)$. If one of the $y_i$,
$i=2,\ldots, k$, is not the identity, we can write $x_1^{-w}y_1$
as the product $x_1^{-w}y_1 = y_i^{v_1}\cdots y_i^{v_e}$ for $e
\leqslant {\rm cn}(G)$, and get the tuple $\mathsf{y}$ after $e$
pairs of moves of the form (below $v_0 = 1$):
\begin{eqnarray*}
(z_1, \ldots, z_k) & \longrightarrow &  (z_1, \ldots,
z_i^{v_{j-1}^{-1}v_j},
 \ldots,
z_k), \; j = 1,2,\ldots, e \\[1ex]
(z_1, \ldots, z_k)  & \longrightarrow & (z_1z_i,z_2, \ldots, z_k)
\end{eqnarray*}
and the correction
$$
(z_1, \ldots, z_k)  \longrightarrow   (z_1, \ldots,
z_i^{v_e^{-1}},  \ldots, z_k).
$$

If, however, all $y_i = 1$, $i =2, \ldots, k$, then $y_1 \ne 1$.
Arguing as before,  $2d+1$ moves will suffice to make $(x_1,
x_1^{-1}y_1, y_3,\ldots,y_k) = (x_1, x_1^{-1}y_1, 1,\ldots,1)$
from $(x_1,1,\ldots, 1)$,  one extra move to make
$(y_1,x_1^{-1}y_1, y_3, \ldots, y_k)$, and at most $2{\rm
cn}(G)+1$ moves to replace $x_1^{-1}y_1$ by $y_2 =1$.

Hence the vertex $\mathsf{x}$ can be connected to $\mathsf{y}$ by
a path of at most
$$
\max\{\,2{\rm or}(G)(k-1)+ 2{\rm cn}(G) +1, \; 2{\rm or}(G)+ 1 +
2{\rm cn}(G)+ 1 \,\} \leqslant 2k{\rm or}(G) + 2{\rm cn}(G)
$$
edges,
 and the proposition follows.
 \hfill $\square$

 \medskip

 To complete the proof of the theorem, we need to list some of the
 known estimates for the covering numbers of finite simple
 groups. They show, in particular, that there is considerable
 scope for improvement of our rather crude estimates.

(a) If $G$ is a Chevalley or twisted Chevalley group of Lie rank
${\rm rank}\,G$ then
$$
{\rm cn}(G) < d \cdot {\rm rank}\,G
$$
for some constant $d$ which does not depend on $G$ ( Ellers,
Gordeev and Herzog \cite{EGH}).

(b) In the case of $PSL_n(q)$, $q\geqslant 4$, $n \geqslant 3$,
there is a better bound
$$
{\rm cn}(PSL_n(q)) = n
$$
(Lev \cite{Le}), while
$$
{\rm cn}(PSL_2(q)) = 3
$$
for all $q \geqslant 4$ (Arad, Chillag and Morgan \cite{ACM}).

(c) For the alternating groups,
\begin{eqnarray*}
{\rm cn}({\rm Alt}_n) & = &  \left[\frac{n}{2}\right], \;\; n
\geqslant
6,\\
{\rm cn}({\rm Alt}_5) & = & 3
\end{eqnarray*}
(Dvir \cite{Dv}).

(d) Covering numbers of sporadic groups do not influence our
asymptotic results. However, it worth mentioning that they are
known (Zisser \cite{Zi}).

 In all these cases ${\rm cn}(G) < c \log_2 |G|$. Since, obviously,
 ${\rm or}(G) \leqslant {\rm cn}(G)$, the theorem follows.
 \hfill $\square$ $\square$

 \bigskip

Modulo the Ore-Thompson Conjecture the estimate of
Proposition~\ref{prop:refined}
 takes the form
$$
{\rm diam}\,(\Delta_k(G,G)) \leqslant 8k + 4{\rm cn}(G).
$$
Notice that the Ore-Thompson Conjecture is true for all Chevalley
groups $G(q)$ for $q \geqslant 8$ (Ellers and Gordeev \cite{EG}).

Also notice that since a simple group $N$ is isomorphically
embedded into its automorphism group,  we have the obvious
inequality
$$
{\rm diam}\,(\Delta_k(G,N)) \leqslant {\rm diam}\,(\Delta_k(N,N))
$$
for every group $G$ which contains $N$ as a normal subgroup.

\medskip

 For the purpose of generating pseudorandom elements in a
simple normal subgroup $N$ of the group $G$, the Andrews-Curtis
graph has the following apparent advantages over the product
replacements graph $\Gamma_k(N)$.  The fact that the diameter
${\rm diam}\,(\Delta_k(G,N))$ is much smaller than that of
$\Gamma_k(N)$ suggests the possibility that the mixing time of a
random walk on $\Delta_k(G,N)$ is smaller than the mixing time of
a random walk on $\Gamma_k(N)$. The vertices of $\Delta_k(G,N)$
are all tuples in $N^k \smallsetminus (1,\ldots, 1)$, hence the
sample of elements of $N$ obtained by taking random components of
random generating tuples (vertices) in $\Delta_k(G,N)$  is not
biased.

\subsection{Gaschuetz's lemma for normal generation}
\label{sec:Gaschuetz}

Notice that an epimorphism $G \longrightarrow H$ of groups induces 
an adjacency-preserving map of graphs $\Delta_k(G,G) 
\longrightarrow \Delta_k(H,H)$. It follows that  the preimage of 
every connected component of $\Delta_k(H,H)$ is the union of some 
connected components of $\Delta_k(G,G)$.

\begin{proposition} {\rm (V.~D.~Mazurov)}
If a finite group\/ $G$ is generated as a normal subgroup by $k$ 
elements\/ {\rm (}that is, $$G = \langle h_1^G,\ldots, 
h_k^G\rangle$$ for some $h_1,\ldots, h_k \in G${\rm )} and the 
images $\bar g_1,\ldots, \bar g_k$ of some elements $g_1,\ldots, 
g_k$ in the factor group $G/M$ for some normal subgroup $M \lhd 
G$ generate $G/M$ as a normal subgroup,
$$G/M =
\langle \bar g_1^{G/M},\ldots, \bar g_k^{G/M}\rangle$$ then there 
exist elements $m_1,\ldots, m_k$ in $M$ such that
$$G = \langle (g_1m_1)^G,\ldots, (g_km_k)^G\rangle.$$
\label{mazurov}
\end{proposition}

\paragraph{Proof.}  The proof is based on the
following classical result:

\begin{fact} {\rm (Gaschuetz \cite{gaschuetz})}
If a finite group\/ $G$ is generated by\/ $k$ elements and the 
images  of some elements\/ $g_1,\ldots, g_k$ in the factor 
group\/ $G/M$ for some normal subgroup\/ $M \lhd G$ generate 
$G/M$, then there exist elements\/ $m_1,\ldots, m_k$ in $M$ such 
that\/ $\langle g_1m_1,\ldots, g_km_k\rangle = G$.
 \label{fact:gaschuetz}
\end{fact}

 In a minimal counter-example to Proposition~\ref{mazurov}, 
 $M$ is a minimal normal subgroup. Let
$H=\langle g_1,\ldots, g_k\rangle_n$ (where $\langle \;\rangle_n$ 
denotes the generation as a normal subgroup in $G$). Then $H \cap 
M = 1$ and $G=H\times M$, so $M$ is simple. If $M$ is non-abelian 
then, for $1\ne m\in M$, $G=\langle g_1m,\ldots, g_k\rangle_n$, 
so $M$ is abelian and  intersects $[G,G]=[H,H]$ trivially. It is 
obvious that $G/[H,H]M=\langle g_1,\ldots, g_k\rangle[H,H]M $.  
By Gasch\"utz lemma (Fact~\ref{fact:gaschuetz}), $G/[H,H]= 
\langle g_1m_1,\ldots, g_km_k\rangle[H,H]$ for some  $m_1,\ldots, 
m_k$ from $M$. These $m_1,\ldots, m_k$ are required elements.
 \hfill $\square$

\medskip
As an immediate corollary we have the following covering property 
of Andrews--Curtis graphs.

\begin{corollary}
 \label{co:epi}
 If\/ $G$ is a finite group normally generated by $k$ elements and $M \lhd G$
  then the canonical map
 $$\Delta_k(G,G) \longrightarrow \Delta_k(G/M,G/M)$$
 is surjective.
\end{corollary}

If the canonical map $\Delta_k(G,G) \longrightarrow 
\Delta_k(H,H)$ is surjective we shall say that the graph 
$\Delta_k(G,G)$ {\em covers} the graph $\Delta_k(H,H)$.

\subsection{The Andrews-Curtis graphs of finite soluble groups}
\label{sec:soluble}

 It is easy to see that the graph
$\Delta_k(G,G)$ is not necessary connected. Indeed, notice that if
$G$ is an abelian group then $\Delta_k(G,G) = {\tilde
\Gamma}_k(G)$. Therefore the following fact is also applicable to
the Andrews-Curtis graphs of abelian groups.

\begin{fact}
Let $A$ be a finite abelian group represented as
$$
A \simeq {\mathbb{Z}}_{e_1} \times\cdots\times  {\mathbb{Z}}_{e_r},
$$
where $e_1 \mid e_2 \mid \ldots \mid e_r$, Then
\begin{itemize}
\item[{\rm (i)}] {\rm (Neumann and Neumann \cite{nn})} $\Gamma_k(A)$ is connected if $k > r$.
\item[{\rm (ii)}] {\rm (Diaconis and Graham \cite{dg})} $\Gamma_r$ has $\varphi(e_1)$
components of equal size.
\end{itemize}
Here  $\varphi(n)$ is the Euler function, i.e.\ the number of
positive integers smaller than $n$ and coprime to $n$.
\label{fact:AC-elementary-abelian}
\end{fact}

However, we shall show in this section that
abelian factor groups is the only obstacle to the
connectedness of the Andrews-Curtis graph of finite soluble
groups.

\begin{lemma}
\label{le:solv} Let\/ $G$ be a soluble {\rm (}not necessary 
finite{\rm )} group. A subset\/ $Y \subset G$ generates $G$ as a 
normal subgroup if and only if\/ $Y$ generates $G$ modulo 
$[G,G]$, i.e., the canonical image of\/ $Y$ generates the 
abelianisation $Ab(G) = G/[G,G]$.
\end{lemma}

\paragraph{Proof} Let $H = \langle Y^G \rangle \lhd G$. Suppose
that $c$ is the derived length of $G$ and $G^{(c)}$ is the last
 (non-trivial) term of the derived series of $G$.  By induction on
 $c$ we may assume that $G = HG^{(c)}$. Now
  $$G/H = HG^{(c)}/H \simeq G^{(c)}/H \cap G^{(c)}$$
  which shows that $G/H$ is abelian. Hence $[G,G] \leq H$.
  Therefore $H = G$, as required. \hfill $\square$ 

\medskip

\begin{corollary}
 If\/ $G$ is a finite soluble group generated as a normal subgroup by $k$ elements
  than the canonical map
 $$\Delta_k(G,G) \longrightarrow \Gamma_k(G/[G,G],G/[G,G])$$
 is surjective.
\end{corollary}

Lemma \ref{le:solv}  allows one to compute the probability
$\psi_k(G)$ that $k$ uniformly and independently chosen elements
in $G$  generate $G$ as a normal subgroup, i.e.,
 $$\psi_k(G) = \frac{|V_k(G,G)|}{|G|^k}.$$
 Observe that if $G$ is a finite abelian group then $\psi_k(G)$
 is just the probability that  $k$ elements from $G$ generate $G$.

\begin{corollary} Let\/ $G$ be a finite soluble group. Then
$$\psi_k(G) = \psi_k(Ab(G)).$$
\end{corollary}

We can now analyse the behaviour of Andrews-Curtis graphs of
finite soluble groups.

\begin{proposition}
Suppose that a finite soluble group $G$ is  generated by $k$
elements $x_1,\ldots ,x_k$. Then for any $f_i\in [G,G]$ the
$k$-tuple $(x_1,\ldots ,x_k)$ is connected by Andrews--Curtis
transformations to $(x_1f_1,\ldots, x_kf_k)$.\/
\label{prop:soluble}

\end{proposition}

\paragraph{Proof.} This is effectively a word-by-word reproduction
of the argument from \cite{My}. We use induction on the derived
length of $G$. Let $A$ be the last non-trivial term of the
derived series of $G$. By the induction hypothesis applied to the
images of the $x_i$ and the $x_if_i$ the corresponding $k$-tuples
in $G/A$ are equivalent. By \cite[Property~1]{My} this implies
that $(x_1f_1,\ldots ,x_kf_k)$ is equivalent (in $G$)   to
$(x_1a_1,\ldots ,x_ka_k)$ for some $a_i\in A$. It remains to
connect $(x_1a_1,\ldots ,x_ka_k)$ with $(x_1,\ldots ,x_k)$ by
Andrews--Curtis transformations. This is done in a process of
successive ``elimination" of the factors $a_i$. At each step the
system of generators $x_i$ of the group $G$ is replaced by
another one by some Nielsen transformations. Let $A_1$ be the
normal closure in $G$ of the elements $a_1,\ldots , a_{k-1}$ and
let the bar denote the images in $G/A_1$. Then $(\bar x_1\bar
a_1,\ldots ,\bar x_k \bar a_k)= (\bar x_1,\ldots ,\bar
x_{k-1},\bar x_k \bar a_k)$. Since the $ \bar x_i$ generate
$\,\overline{\! G}$, by \cite[Property~2]{My} $(\bar x_1,\ldots
,\bar x_{k-1}, \bar x_k \bar a_k)$ is connected to $(\bar
x_1,\ldots ,\bar x_{k-1}, \bar x_k)$ by a chain of
Andrews--Curtis transformations applied only to the last
component. Lifting these transformations to $G$ we obtain that
$(x_1a_1,\ldots ,x_ka_k)$ is equivalent to $(x_1a_1,\ldots
,x_{k-1}a_{k-1},x_ka'_k)$, where $a_k'\in A_1$. We write
$a_k'=a_{j_1}^{g_1}\cdots a_{j_l}^{g_l}$, where $ 1\leqslant
j_s\leqslant k-1$ and the $g_{s}$  are some elements of $G$. Then
we successively kill all the factors $ a_{j_s}^{g_s}$ in the last
component at the expense of changing $x_n$ by some Nielsen
transformations. Namely, by \cite[Property~4]{My} $(x_1a_1,\ldots
,x_ka'_k)$ is equivalent to $(x_1a_1,\ldots
,x_{j_l}a_{j_l}^{g_l},\ldots ,x_ka'_k)$. Then we apply
Andrews--Curtis transformations to the last component to get
$$(x_1a_1,\ldots ,x_{j_l}a_{j_l}^{g_l},\ldots ,x_ka'_k
(a_{j_l}^{-1})^{g_l} x_{j_l}^{-1})= (x_1a_1,\ldots
,x_kx_{j_l}^{-1} a_{j_1}^{g'_1}\cdots
a_{j_{l-1}}^{g'_{{l-1}}}),$$ where $
g'_{s}=g_{s}^{x_{j_l}^{-1}}$. The number of ``$a^g$-factors" in
the last component is now smaller, although the generator $x_k$
is ``replaced" by $ x_kx_{j_k}^{-1}$. After finitely many such
steps we arrive at $(x_1a_1,\ldots ,x_{k-1}a_{k-1},x_k')$, where
$x_k'$ is a result of Nielsen transformations. We have thus got
rid of the $a$-factor in the last component. This process can now
be repeated with $A_1$ generated by fewer elements, although  for
a new system of generators of $G$ obtained from the $x_i$ by
Nielsen transformations. It is of course important that at each
step we are dealing with a $k$-tuple of the form $(x'_1a_1,\ldots
,x'_{k}a_{k})$, where the $x_i'$ are generators of $G$. To use
formally an induction argument, one can each time simply
rearrange the components of the $k$-tuple, which can be done by
the Andres--Curtis transformations, so that some initial segment
of components increasing in length is
  free of
$a$-factors. The last step of this process is also covered by
this argument, when \cite[Property~2]{My} is applied as above
with $A_1=1$.   Finally  we shall arrive at an equivalent
$k$-tuple $(x_1',\ldots ,x'_k)$ obtained from $ (x_1,\ldots
,x_k)$  by Nielsen transformations. Reversing the chain of these
Nielsen transformations we arrive at $(x_1,\ldots ,x_k)$. Thus,
$(x_1a_1,\ldots ,x_ka_k)$ and therefore $(x_1f_1,\ldots ,x_kf_k)$
is equivalent to $(x_1,\ldots ,x_k)$. The proposition is proved.
\hfill $\square$

\medskip

In Proposition~\ref{prop:soluble}, it would be interesting to 
replace `generated by $k$ elements' by `generated as a normal 
subgroup by $k$ elements':

\begin{question}
\textsf{Suppose that a finite soluble group $G$ is  generated as 
a normal subgroup by $k$ elements $x_1,\ldots ,x_k$. Is it true 
that for any $f_i\in [G,G]$ the $k$-tuple $(x_1,\ldots ,x_k)$ is 
connected by Andrews--Curtis transformations to $(x_1f_1,\ldots, 
x_kf_k)$?\/ \label{q:soluble} }
\end{question}

\begin{theorem}
Suppose that  a finite soluble group $G$  can be generated by $k$
elements. Then the preimages in $G$ of the connected components
of the Andrews--Curtis $k$-tuple graph of $G/[G,G]$ are all
connected. In particular, the graph $\Delta_{k+1}(G,G)$ is
connected.   \label{th:soluble}
\end{theorem}

\paragraph{Proof.} Let $(u_1,\ldots ,u_k)$ and $(v_1,\ldots ,v_k)$ be two $
k$-tuples of elements of $G$ (each generating $G$ modulo $[G,G]$)
that are equivalent modulo $[G,G]$. Then $(u_1,\ldots ,u_k)$ is
equivalent in $G$ to $(v_1f_1,\ldots ,v_kf_k)$ for some $f_i\in
[G,G]$.  Since $G$ is $k$-generated, by Gasch\"utz' lemma there
are elements $ h_i\in [G,G]$ such that the elements
$v_1h_1,\ldots ,v_kh_k$ generate $ G$. By
Proposition~\ref{prop:soluble} the $k$-tuples $(v_1,\ldots
,v_k)$ and $(v_1h_1,\ldots ,v_kh_k)$ are equivalent, as well as
the $k$-tuples $(v_1h_1,\ldots ,v_kh_k)$ and $(v_1f_1,\ldots
,v_kf_k)$. By transitivity hence $(u_1,\ldots ,u_k)$ and
$(v_1,\ldots ,v_k)$ are equivalent, as required. Since the factor
group $G/[G,G]$ is $k$-generated, the graph
$\Delta_{k+1}(G/[G,G],G/[G,G])$ is connected by
Fact~\ref{fact:AC-elementary-abelian}, and the second assertion
also follows. \hfill $\square$

\medskip

One can compare this result with the following observation about
soluble groups.

\begin{fact} {\rm (Dunwoody \cite{dunwoody}, see also \cite[Theorem~2.36]{P2})}
Let $G$ be a finite soluble group generated by $k$ elements.
Then $\Gamma_{k+1}(G)$ is connected.

\end{fact}

Theorem~\ref{th:soluble} is especially nice in the case of
nilpotent groups. A well-known fact about nilpotent groups
states that if  $A$ is a subgroup of a nilpotent group $G$ such
that $A[G,G]=G$, then $A=G$. It follows that a tuple of elements
generates $G$ as a normal subgroup if and only if it generates it
as a group. Applying Theorem~\ref{th:soluble}, we immediately
have

\begin{theorem}
Suppose that $G$ is a finite nilpotent group. Then the pre-images
in $G$ of the connected components of the Andrews--Curtis graph
of\/ $G/[G,G]$ are all connected.

 \label{th:nilpotent}
\end{theorem}

\paragraph{Remarks.} In fact, an argument in \cite{My} states that the tuples
$(g_1,\ldots ,g_k)$ and  $(g_1,\ldots
,g_{k-1},g_kf)$ for $f\in [G,G]$ are connected by Andrews--Curtis
transformations in any group $ G=\left< g_1,\ldots ,g_k\right>$. In
the above case of $G$ being nilpotent there is also an alternative
computation based on using more of commutator calculus, which
may
be more efficient from computational viewpoint.

\section{The Andrews-Curtis Conjecture:
an approach via unsoluble finite groups}
\label{sec:attack}
\subsection{Disconnected Andrews-Curtis graphs}

A possible way to confirm a counterexample to the Andrews-Curtis
Conjecture starts with  one of the suggested potential
counterexamples, that is, words $u=u(x,y)$ and $v=v(x,y)$ which
generate the free group $F_2 = \langle x,y\rangle$ of rank two,
and which are suspected of not being connected to $x$ and $y$ by
a sequence of Andrews-Curtis moves. One can take a finite group
$G$ with more than one connected component of $\Delta_2(G,G)$ and
consider the map
\begin{eqnarray*}
\omega: \Delta_2(G,G) & \longrightarrow & \Delta_2(G,G)\\
        (x,y) & \mapsto & (u(x,y),v(x,y)).
\end{eqnarray*}
If $\omega$ maps a vertex $(x,y)$ to a vertex which belongs to a
different  component of $\Delta_2(G,G)$, then the pairs $(x,y)$
and $(u(x,y),v(x,y))$ obviously constitute a counterexample to the
Andrews-Curtis Conjecture. Of course, the first candidate for the
map $\omega$ should come from the simplest possible potential
counterexample to the AC-conjecture. The following pair
\begin{equation}
\label{eq:k}
 (x^3y^{-4}, xyxy^{-1}x^{-1}y^{-1})
 \end{equation}
 occurs in the second presentation in  the series  of
potential counterexamples proposed by Akbulut and Kirby \cite{AK}.
The total length of these words is equal to 13.  Note that all
pairs $(u,v)$ which generate $F_2$ as a normal subgroup and have
the total length $|u| + |v| \leq 12$ satisfy the AC-conjecture
\cite{My2}. So the potential counterexample (\ref{eq:k})  has the
minimal possible length. Moreover, recently, Havas and Ramsay
proved that every pair of elements in $F_2$ which generates $F_2$
as a normal subgroup and has the total length 13 is AC-equivalent
either to $(x,y)$ or to (\ref{eq:k}) \cite{Havas}. This shows that
the map
$$
\omega: (x,y) \mapsto (x^3y^{-4}, xyxy^{-1}x^{-1}y^{-1})
$$
should be of prime interest here.

However, in view of Myasnikov's result on the Andrews-Curtis
graphs of free soluble groups (Fact~\ref{fact:myasnikov}), we
should not expect to find a counterexample to the Andrews-Curtis
Conjecture among finite soluble groups.

For that reason it would be interesting to study the
Andrews-Curtis graph for finite {\em unsoluble} groups.

\paragraph{Perfect groups.} 
Recall that a group $G$ is called {\em perfect} if it coincides
with its commutator, $G = [G,G]$.

\begin{question}
 \textsf{ Is it true that, for a perfect finite group $G$,
the Andrews--Curtis graph  $\Delta=\Delta_k(G,G)$ is connected?}

\textsf{More generally, is it true that the preimage in
$\Delta_k(G,G)$ of every connected component of
$\Delta_k(G/[G,G], G/[G,G])$ is connected? } \label{q:connected}
\label{q:perfect=connected}
\end{question}

In the case of a simple finite group $G$ the answer is obviously
`yes'. Moreover, a slightly more general result is true.

\begin{lemma} {\rm (Sukru Yalcinakya)}
If a group $G$ has a unique maximal normal subgroup $M$ then 
$\Delta_k(G,G)$ is connected for every $k \geqslant 2$. 
\label{lm:sukru}
\end{lemma}

\paragraph{Proof.} Notice that every element $x \in G \smallsetminus M$ 
generates $G$ as a normal subgroup, and vice versa. Therefore the 
vertices of $\Delta_k(G,G)$ are all $k$-tuples $(x_1,\ldots, 
x_k)$ such that at least one of $x_i$ does not belong to $M$. Now 
any two vertices can be obviously connected by the 
Andrews--Curtis moves. \hfill $\square$

\medskip

The positive answer to Question~\ref{q:connected} would 
implicitly suggest that the Andrews-Curtis Conjecture is true. 
However, the connectedness of the Andrews--Curtis graphs is not 
yet the end of the story, see Question~\ref{q:distance-grows} 
below.

\begin{question}
\textsf{Let\/ $K$ and  $L$ be normal subgroups of a finite 
grou\/p  $G$ and\/ $K \cap L =1$. Assume that the Andrews--Curtis 
graphs\/ $\Delta_k(G/K,G/K)$ and\/ $\Delta_l(G/L,G/L)$ are 
connected. Is it true that\/ $\Delta_{k+l}(G,G)$ is connected? } 

\textsf{In more general terms, how do connected components of\/  
$\Delta_{m}(G,G)$ relate to  connected components of\/ 
$\Delta_k(G/K,G/K)$ and\/ $\Delta_l(G/L,G/L)$?}
\end{question}

The first part of the question is likely to be easy.

\paragraph{Non-perfect groups.}
If the preimage $\Gamma$ of a connected component from
$\Delta_2(G/[G,G],G/[G,G]))$ is disconnected, it would be very
interesting to look at the map
$$
\omega: (x,y) \mapsto (u(x,y),v(x,y))
$$
on $\Gamma$. Notice that since the images $\bar u$ and $\bar v$ of
$u(x,y)$ and $v(x,y)$ generate the  factor group
$A=F_2/[F_2,F_2]$ of the free group $F_2 = \langle x,y\rangle$
modulo the commutator, the determinant $\det(\tilde u,\tilde v)$
of the matrix representing $\tilde u, \tilde v$ in the basis
$\tilde x, \tilde y$ is $\pm 1$. Hence the image $(\bar u, \bar
v)$ of $\omega(x,y)$ belongs to the same connected component of
$\Delta_2(G/[G,G],G/[G,G])$ as $(\bar x, \bar y)$, and
$\omega(\Gamma) \subseteq \Gamma$.

\begin{question}
\textsf{Does the map $\omega$ preserve the connected components
of $\Gamma$?}
\end{question}

The negative answer to this question, of course, provides a
counterexample to the Andrews--Curtis Conjecture. So far this
gives only a one way approach to the problem. It would be
interesting to see if there exists a two-way reduction of the
Andrews--Curtis Conjecture to questions about Andrews--Curtis
graphs of finite groups. In particular, we would like to mention
the following conjecture.

\begin{conjecture}
\textsf{ If the normal generators\/ $x, y$ and\/ $u(x,y), v(x,y)$
of\/ $F_2$ are not connected by Andrews--Curtis transformations,
then there exists a finite factor group $G$ where $(\bar x, \bar
y)$ and $(u(\bar x, \bar y), v(\bar x, \bar y)$ belong to
different connected components of $\Delta_2(G,G)$.}

\end{conjecture}

\subsection{\relax An alternative approach to construction of\\ counterexamples}

As we have already mentioned, good connection  properties of the
Andrews--Curtis graphs of finite groups do not yet herald the end
of attempts to construct a counterexample to the Andrews--Curtis
Conjecture by means of finite group theory. We can also try an
alternative approach.

Notice that if the normal generators $(x,y)$ and $u(x,y), v(x,y)$
of\/ $F_2$ are connected by $d$ Andrews--Curtis transformations,
then the same is true for an arbitrary finite group. Therefore we
come to the following question.

\begin{question}
\textsf{Does there exists a series of finite groups\/ $\{G_n\}$
such that the path distance $d_n$ in $\Delta_2(G_n,G_n)$ between
the pairs\/ $(\bar x, \bar y)$ and\/ $(u(\bar x, \bar y), v(\bar
x, \bar y))$  is unbounded?} \label{q:unbounded}
\end{question}

\begin{question}
\textsf{Assuming that the Andrews--Curtis Conjecture is false 
and  the normal generators\/ $x, y$ and $u(x,y), v(x,y)$ of\/ 
$F_2$ are not connected by Andrews--Curtis transformations, is it 
true  there exists a series of finite factor groups\/ $\{G_n\}$ 
such that if\/ $d_n$ is  the path distance in $\Delta_2(G_n,G_n)$ 
between the pairs\/ $(\bar x, \bar y)$ and\/ $(u(\bar x, \bar y), 
v(\bar x, \bar y))$ then the sequence $\{d_n\}$ is unbounded? } 
\label{q:distance-grows}
\end{question}

Of course, similar conjectures can be formulated for arbitrary
lengths $k$ of $k$-tuples of normal generators.

Notice, however, that the free group $F_n$ is residually $\cal S$ 
for every infinite set $\cal S$ of pairwise nonisomorphic finite 
nonabelian simple groups \cite{weigel1,weigel2,weigel3}, that is, 
for each $g\in F_n$, $g\not= 1$, there is an epimorphism 
(depending on $g$) from $F_n$ onto a group in $\cal S$ such that 
the image of $g$ is not $1$. This means, for example, that $F_n$ 
is residually ${\rm Chev}_n$ for the class ${\rm Chev}_n$ of all 
finite Chevalley groups of rank $\leqslant n$, while these groups 
have uniformly bounded diameters of their Andrews--Curtis graphs 
(Section~\ref{sec:3}). 

Therefore Questions~\ref{q:unbounded} and \ref{q:distance-grows} 
might happen to be hard to resolve. However, their versions for 
{\em restricted Andrews--Curtis} graphs are more likely to have 
positive solutions and seem to be more accessible for a study by 
means of computer experiments.

\subsection{Restricted Andrews-Curtis graphs}

Let $G$ be a group generated by a set $S$ and  $N \lhd G$. We 
define the {\em restricted Andrews--Curtis graph\/} 
$\bar\Delta_k(G,S,N)$  as the graph with the same vertices as in 
$\Delta_k(G,N)$, that is, $k$-tuples of elements in $N$ which 
generate $N$ as a normal subgroup:
$$ \left\{\, (h_1,\ldots, h_k) \mid \langle
h_1^G,\ldots, h_k^G\rangle = N \,\right\}.$$ 
Two vertices are 
connected by an edge if one of them is obtained from another by 
one of the moves:
\begin{eqnarray*}
(x_1, \ldots, x_i, \ldots, x_j, \ldots, x_k) & \longrightarrow & (x_1, \ldots,x_ix_j^{\pm 1}, \ldots,  x_k), \; i \ne j \\[2ex]
(x_1, \ldots, x_i, \ldots, x_j, \ldots, x_k) & \longrightarrow & (x_1, \ldots,x_j^{\pm 1}x_i, \ldots,  x_k), \; i \ne j \\[2ex]
(x_1, \ldots, x_i, \ldots, x_k) & \longrightarrow &  (x_1,\ldots,x_i^{s}, \ldots, x_k), \;  s\in S\\[2ex]
(x_1, \ldots, x_i, \ldots, x_k) & \longrightarrow & (x_1, 
\ldots,x_i^{-1}, \ldots,  x_k).
\end{eqnarray*}

Thus $\bar\Delta_k(G,S,N)$ is a subgraph of $\Delta_k(G,N)$ whose 
edges correspond to Nielsen moves, inversions, and to conjugation 
by generators $s\in S$ rather than arbitrary elements $w\in G$.

Obviously, the graphs $\bar\Delta_k(G,S,N)$ and $\Delta_k(G,N)$ 
have the same connected components. For finite groups $G$, the 
graph $\bar\Delta_k(G,S,N)$ has much large diameter than  
$\Delta_k(G,N)$. Of course, we have the obvious estimate 
$$
{\rm diam}\,  \bar\Delta_k(G,S,N) \leqslant {\rm 
diam}\,\Delta_k(G,N) \cdot {\rm diam}\, \mathsf{Cayley}(G,S),
$$
where $\mathsf{Cayley}(G,S)$ is the Cayley graph of the group $G$ 
with respect to the generating set $S$. 

\begin{question}
\textsf{Find better bounds for the diameter of the restricted 
Andrews--Curtis graph\/ $\bar\Delta_k(G,S,N)$ in the case of 
finite simple groups $G=N$ and `natural' sets of generators.}
\end{question}

Notice that if the normal generators $(x,y)$ and $u(x,y), v(x,y)$ 
of\/ $F_2$ are connected by $d$ edges in the restricted 
Andrews--Curtis graph of $F_2$ with respect to some generating 
set $\{a,b\}$ of $F_2$, then the same is true for an arbitrary 
finite group. Therefore we come to the following analogues of 
Questions~\ref{q:unbounded} and \ref{q:distance-grows}.

\begin{question}
\textsf{Does there exists a series of finite groups\/ $\{G_n\}$ 
with generators $x_n, y_n$ such that the path distance $d_n$ in 
$\bar\Delta_2(G_n,\{x_n,y_n\},G_n)$ between the pairs\/ $( x_n, \ 
y_n)$ and\/ $(u( x_n, y_n), v( x_n,  y_n))$  is unbounded?} 
\label{qq:unbounded}
\end{question}

\begin{question}
\textsf{Assume that the Andrews--Curtis Conjecture is false and  
the normal generators\/ $u(x,y), v(x,y)$  of the free group $F_2$ 
are not connected by Andrews--Curtis transformations to the free  
generators $x$ and $y$. Is it true that  there exists a series of 
finite factor groups\/ $\{G_n\}$   $F_2$  such that if\/ $d_n$ 
is  the path distance in $\bar\Delta_2(G_n,\{\bar x, \bar 
y\},G_n)$ between the pairs\/ $(\bar x, \bar y)$ and\/ $(u(\bar 
x, \bar y), v( \bar x,  \bar y))$ then the sequence $\{d_n\}$ is 
unbounded? } \label{qq:distance-grows}
\end{question}

It would be interesting to try to run computer experiments with
the restricted  Andrews--Curtis graphs of finite simple groups as 
an attempt to analise their metric properties. Taking, for 
example, the transvections
$$
x = \left(\begin{array}{cc} 1 & 0 \cr 2 & 1\end{array}\right)
\;\hbox{ and }\; y = \left(\begin{array}{cc} 1 & 2 \cr 0 &
1\end{array}\right) $$ in the group $G_q = SL_2({\mathbb{F}}_q)$
for reasonably small values of $q$, is it feasible to compute the
path distance $d_q$ in $\bar\Delta_2(G_q,G_q)$ from $(x,y)$ to
$(x^3y^{-4}, xyxy^{-1}x^{-1}y^{-1})$? Might it happen that  a
geodesic path found in $SL_2({\mathbb{F}}_q)$ can be lifted to 
the free group 
 $$
\left\langle \left(\begin{array}{cc} 1 & 0 \cr 2 & 
1\end{array}\right),   \left(\begin{array}{cc} 1 & 2 \cr 0 & 
1\end{array}\right)\right\rangle \leqslant SL_2({\mathbb{Z}})? $$

 However, even if the growth of $d_q$ is
detected in a small sample of computationally accessible graphs,
we still encounter a possibly very difficult problem of
theoretical analysis of metric properties of (restricted) 
Andrews--Curtis graphs of arbitrary big size.

\subsection{Expanders}

The following result by Lubotzky and Pak is intimately related to 
their\linebreak Theorem~\ref{th:LuP}.

\begin{theorem} {\rm \cite{LuP}}
If\/ $Aut \, F_k$ has property {\rm (}T{\rm )} and $G$ is a 
finite group then\/ every connected component of\/  $\Gamma_k(G)$ 
is an $\varepsilon$-expander for some $\varepsilon$ which depends 
only on $k$.
\end{theorem}

Here, a graph $\Gamma$  is an {\em $\varepsilon$-expander\/} if, 
for every set of vertices $B \subset \Gamma$ which is less than 
half of $\Gamma$, $|B| < \frac{1}{2} |\Gamma|$, has sufficiently 
many `new' neighbours:
$$
\left|\left\{\begin{array}{c} \hbox{vertices connected}\\
\hbox{to } B, \hbox{ but not in } B
\end{array} \right\}\right| \geqslant \varepsilon \cdot |B|.
$$

We would like to mention, in passing, the following question.

\begin{question}
\textsf{Do the connected components of the graphs\/ 
$\bar\Delta(G,S,N)$ for normal subgroups $N$ of finite groups $G$ 
form a family of expanders?}

\end{question}

\begin{center}
{\bf Acknowledgement}
\end{center}

The authors wish to thank Robert Gilman who was at the origins of 
this project and contributed many important ideas, and Victor D. 
Mazurov, who proved Proposition~\ref{mazurov}. We also thank 
Sukru Yalcinkaya for an useful observation which became 
Lemma~\ref{lm:sukru}.

\bigskip

\vfill

 \noindent \textsf{Alexandre V. Borovik, Department of
Mathematics, UMIST, PO Box 88, Manchester M60 1QD, United Kingdom}

\noindent {\tt borovik@umist.ac.uk}

\noindent {\tt http://www.ma.umist.ac.uk/avb/}

\medskip
\noindent \textsf{Evgenii I. Khukhro,
Institute of Mathematics, Novosibirsk-90, 630090, Russia}

\noindent {\tt khukhro@cardiff.ac.uk}

\medskip

\noindent \textsf{Alexei G. Myasnikov,
 Department of Mathematics, The City  College  of New York, New York,
NY 10031, USA}

\noindent  {\tt alexeim@att.net}

\noindent {\tt http://home.att.net/\~\,alexeim/index.htm}
\end{document}